\documentclass{amsart}
\usepackage{amsfonts,amscd}

\usepackage[arrow,matrix,curve,cmtip,ps]{xy}
\newtheorem{theorem}{Theorem}[section]
\newtheorem{lemma}[theorem]{Lemma}
\newtheorem{corollary}[theorem]{Corollary}
\newtheorem{proposition}[theorem]{Proposition}
\theoremstyle{remark}

\theoremstyle{definition}
\newtheorem{definition}[theorem]{Definition}

\numberwithin{equation}{section} \makeatother
\DeclareMathOperator{\Bdb}{{\mathbb B}}

\DeclareMathOperator{\Cdb}{{\mathbb C}}

\DeclareMathOperator{\Ndb}{{\mathbb N}}

\begin{document}

\title{A characterization and a generalization of $W^*$-modules}

\author[David P. Blecher]{David P. Blecher*}
\address{Department of Mathematics, University of Houston, Houston, TX
77204-3008} \email{dblecher@math.uh.edu}

\author[Upasana Kashyap]{Upasana Kashyap}
\address{Department of Mathematics, University of Houston,
Houston, TX  77204-3008} \email{upasana@math.uh.edu}
\date{\today}
\thanks{*Supported by grant DMS 0400731 from the National Science Foundation}


\begin{abstract}
We give a new Banach module characterization of $W^*$-modules, also
known as selfdual Hilbert $C^*$-modules over a von Neumann algebra.
This leads to a generalization of the notion, and the theory, of
$W^*$-modules, to the setting where the operator algebras are
$\sigma$-weakly closed algebras of operators on a Hilbert space.
That is, we find the appropriate weak* topology variant of our
earlier notion of {\em rigged modules}, and their theory, which in
turn generalizes the notions of  $C^*$-module, and Hilbert space,
successively.  Our {\em $w^*$-rigged modules} have canonical `envelopes'
which are $W^*$-modules.  Indeed, $w^*$-rigged modules may be
defined to be a subspace of a $W^*$-module possessing certain
properties.
\end{abstract}

\maketitle

\section{Introduction and Notation}

A {\em $W^*$-module} is a Hilbert $C^*$-module over a von Neumann
algebra which is {\em selfdual}, or, equivalently, which has a
predual (see e.g.\ \cite{Zet,EOR,BM}). These objects were first
studied by Paschke, and then by Rieffel \cite{Pas,Rief} (see also
e.g.\ \cite[Section 8.7]{BLM} for an account of their theory). They
are by now a fundamental object in $C^*$-algebra theory and
noncommutative geometry, being intimately related to Connes'
correspondences for example (see e.g.\ \cite{BDH} for the
relationship). $W^*$-modules have many characterizations; the one
mentioned in our title characterizes them in the setting of Banach
modules in a new way. This in turn leads us to generalize the notion
of `$W^*$-module' to the setting of modules over a dual operator
algebra (by which we mean a $\sigma$-weakly closed algebra of
operators on a Hilbert space). The new class of modules 
we call {\em $w^*$-rigged modules}.
These are the appropriate
weak* topology variant  of our earlier notion of {\em rigged
modules} (see e.g.\ \cite{DB2,BMP,DBsur,BHN}), which in turn
generalizes the notions of Hilbert $C^*$-modules and Hilbert space,
successively. 
$W^*$-rigged modules will have an extensive theory, we
will only present the basics here.
Our main motivation for this project is that many ideas in the
$W^*$-module theory, such as `induced representations', are
beautiful and fundamental, and thus should be important in a larger
context. We also wished to enlarge the universe in which
$W^*$-modules reside and act.    In addition, one obtains many
theorems reprising basic facts from the theory of rings and modules,
but which only make sense for modules satisfying our definition (for
example because direct sums are problematic for general operator
modules).

Our theory utilizes several pretty ideas from operator space theory.
For example, our theory of the space of left multipliers ${\mathcal
M}_{\ell}(X)$ of an operator space $X$ (see e.g.\ \cite[Chapter
4]{BLM}), plays a considerable role in our paper. Indeed the absence
of this tool, and also of a recently introduced module tensor
product \cite{EP} (see also \cite[Section 2]{BK1}), is the main
reason why headway was not made on this project many years ago.

Unlike the $W^*$-module situation, $w^*$-rigged modules do not
necessarily give rise to a weak* Morita equivalence in the sense of
our earlier paper \cite{BK1}.  Nor are they all complemented 
submodules of a direct sum of copies of $M$, as is the case 
for  $W^*$-modules \cite{Pas}.
However, each $w^*$-rigged
module has a canonical $W^*$-module envelope, called the
$W^*$-dilation (and thus  $w^*$-rigged modules
give new examples of $W^*$-modules).  This dilation
 is an important tool in our theory. Indeed, a $w^*$-rigged
module may be defined to be a subspace of a $W^*$-module possessing
certain properties.  Thus there is a von Neumann algebra valued
inner product (and Morita equivalence) around.

The material on $W^*$-modules in  Section 2 is closely related to a paper
of the first author from about ten years ago \cite{Bsd}. The
main point of the latter paper
was that $W^*$-modules fall comfortably into a dual operator module
setting; for example, their usual tensor product (sometimes called
`composition of $W^*$-correspondences'), agrees with a certain
operator space tensor product studied by Magajna.  This had certain
advantages, for example new results about this tensor product (see
also \cite{DH}). Here we show for example that this tensor product
also equals the normal module Haagerup tensor product recently
introduced in
 \cite{EP}, and studied further in \cite{BK1}.  In Section 3 we
 find the variant for $w^*$-rigged modules of the basic theory
 of rigged modules from \cite{DB2} (and to a lesser extent, \cite{BMP}).
 In Section 4 we give several alternative equivalent definitions
 of $w^*$-rigged modules, some of which the reader may prefer.
 In Section 5 we give many examples of $w^*$-rigged modules.

Since a number of the ideas and proofs here are quite analogous to
those from our papers on related topics, principally
\cite{DB2,BMP,BK1,BM}, we will be quite brief in many of the proofs.
We assume that the reader has some familiarity with these
earlier ideas and proof techniques, and will often merely indicate
the new techniques which allow one to modify  arguments from
\cite{DB2,BMP} to the present (weak* topology) setting. As the paper
proceeds we will include fewer and fewer details, since we will be
assuming that the reader is growing in familiarity with the methods
introduced earlier and in \cite{BK1}, and the ways these are used to
modify proofs from the older theory.  Also, a few complementary
facts may be found in \cite{UK1,UK}.  Indeed the present paper
represents a research program that the first author suggested to the
second while she was a graduate student,
and it was intended that full
details and other aspects of the theory of $w^*$-rigged modules
not touched upon here, would be presented elsewhere. 
This project is still in progress.

We will use the notation of our previous paper \cite{BK1}. We will
assume that the reader is familiar with basic notions from operator space theory, 
as may be found in any of the current texts
on that subject,  and the application of this theory to operator
algebras (see e.g.\ \cite{BLM}).  The latter source
may also be consulted for background and for any
unexplained terms below.

We also assume that the reader is familiar with basic Banach space
(and operator space) duality principles (see e.g.\ \cite[Section
1.4, 1.6, Appendix A.2]{BLM}). We will often abbreviate `weak*' to
`$w^*$'. Throughout the paper, $H$ and $K$ denote Hilbert spaces.
Unless indicated otherwise, $M$ denotes a {\em dual operator
algebra}, that is, an operator algebra which is also a dual operator
space. We take all dual operator algebras $M$ to be {\em unital},
that is we assume they possess an identity of norm $1$ (this is for
convenience, for nonunital algebras $M$ one may consider the
unitization $M^1$ (see \cite[2.7.4 (5)]{BLM}); for example defining
a module over a nonunital algebra $M$ to be $w^*$-rigged if it is
$w^*$-rigged over $M^1$). By well
 known duality principles, any $w^*$-closed subalgebra of $B(H)$, is a
 dual operator algebra.  Conversely, for any dual
 operator algebra $M$,  there exists a Hilbert
 space $H$ and a $w^*$-continuous completely isometric homomorphism
 $\pi: M \rightarrow B(H)$ (see e.g.\
 \cite[Theorem 2.7.9]{BLM}). Then the range $\pi(M)$ is a
 $w^*$-closed subalgebra of $B(H)$, which we may identify with $M$
in every way.

 For cardinals or sets $I, J$, we use the symbol $M_{I,J}(X)$ for the
 operator space of $I \times J$
 matrices over an operator space $X$, whose `finite submatrices' have uniformly
 bounded norm.    We set $C^w_J(X) = M_{J,1}(X)$ and $R^w_J(X) =
 M_{1,J}(X)$; and these are written as $C_n(X)$ and $R_n(X)$ if $J = n$ is finite.

A  {\em concrete left operator module} over an operator algebra $A$,
 is a linear subspace  $X\subset B(K,H)$, such that $ \pi(A)X$ $\subset X$ for a
 completely contractive representation $ \pi : A \rightarrow B(H)$. An
 {\em abstract operator  A-module} is an operator space $X$ which is also
 an $A$-module, such that $X$ is completely isometrically isomorphic,
 via an $A$-module map, to a concrete operator $A$-module.  Similarly for right
 modules, or
 bimodules.  Most of the
 interesting modules over operator algebras are operator modules, such
 as Hilbert $C^*$-modules (the operator space structure
 on a $C^*$-module is the one it receives as a subspace of its {\em linking
 $C^*$-algebra}---see e.g.\ \cite[Section 8.2]{BLM}).

A {\em normal Hilbert module} over a dual operator algebra $M$ is a
pair $(H, \pi)$, where $H$ is a (column) Hilbert space (see e.g.\
1.2.23 in \cite{BLM}), and $\pi : M \to B(H)$ is a
$w^{*}$-continuous unital completely contractive representation. We
shall call such $\pi$ a {\em normal representation} of $M$. The
module action is given by $m \cdot \zeta = \pi(m) \zeta$.

A {\em concrete dual operator $M$-$N$-bimodule} is a $w^{*}$-closed
 subspace $X$ of $B(K,H)$ such that $\theta(M) X \pi(N)$ $\subset X$, where
 $\theta$ and $\pi$ are normal
 representations of $M$ and $N$ on $H$ and $K$ respectively. An
 {\em abstract
 dual operator $M$-$N$-bimodule} is defined to be a non-degenerate
  operator
 $M$-$N$-bimodule $X$, which is also a dual operator space, such that
 the module actions are separately weak* continuous.  Such
 spaces can be represented completely isometrically as concrete
 dual operator bimodules, and in fact this can be done
 under even weaker hypotheses
 (see e.g.\ \cite{BLM,BM,ER}).  Similarly for
one-sided modules (the case $M$ or $N$ equals $\Cdb$). We use
standard notation for module mapping spaces, e.g.\ $CB(X,N)_{N}$
(resp.\ $w^*CB(X,N)_N$) are the completely bounded (resp.\ and
$w^*$-continuous) right $N$-module maps from $X$ to $N$. Any normal
Hilbert $M$-module $H$ (with its column Hilbert space structure
$H^c$) is a left dual operator $M$-module.  In a couple of proofs,
we will assume that the reader is familiar with the theory of
multipliers of an operator space $X$ (see e.g.\ \cite[Chapter
4]{BLM}). We recall that the left multiplier algebra ${\mathcal
M}_{\ell}(X)$ of $X$ is a collection of certain operators on $X$,
which are weak* continuous if $X$ is a dual operator space
\cite{BM}.  Indeed, in the latter case, ${\mathcal
M}_{\ell}(X)$ is a dual operator algebra.

A bilinear map $u : X \times Y \to Z$ is {\em $M$-balanced} if
$u(xm,y) = u(x, my)$ for $m \in M$, and  {\em completely
contractive} if
it corresponds to a linear complete
contraction on $X \otimes_{h} Y$.   We use the {\em normal module
Haagerup tensor product} $\otimes^{\sigma h}_M$ throughout the
paper, and  its universal property from \cite{EP}, which loosely
says that it `linearizes completely contractive $M$-balanced
separately weak* continuous bilinear maps':
Every completely contractive separately weak* continuous $M$-balanced map $u :
X \times Y \to Z$, induces a
completely contractive weak* continuous complete contraction $X
\otimes^{\sigma h}_{M} Y \to Z$.  We also assume that
the reader is
familiar with notation and facts
about $\otimes^{\sigma h}_M$ from \cite[Section 2]{BK1}.

\section{$W^*$-modules}

We begin this section with a useful lemma:

\begin{lemma}\label{A}
Let $\{ H_{\alpha} \}$ be a collection of Hilbert spaces (resp.\
column Hilbert spaces) indexed by a directed set. Let $Y$ be a dual
Banach space (resp.\ dual operator space). Suppose there exist
$w^{*}$-continuous contractive (resp. completely contractive) linear
maps $\phi_{\alpha} : Y \to H_{\alpha}$, $\psi_{\alpha} : H_{\alpha}
\to Y$, such that $\psi_{\alpha}( \phi_{\alpha}(y)) \buildrel w^{*}
\over \to y$ for each $y \in Y$. Then $Y$ is a Hilbert space (resp.
column Hilbert  space). The inner product on $Y$ is $\langle
y,z\rangle $= {\rm lim}$_{\alpha} \; \langle \phi_{\alpha}(y),
\phi_{\alpha}(z)\rangle$, for $y, z \in Y$.
\end{lemma}

\begin{proof}
The proof that $Y$ is a Hilbert space (resp. column Hilbert space)
follows by the  ultraproduct argument in Theorem 3.10 in \cite{BK1}.
For the last assertion, we will show first that $\lVert
\phi_{\alpha}(y)\rVert^2$ $\to$ $\lVert y\rVert^2$.  Then by the
polarization identity, it follows that $\langle y,z\rangle $=
lim$_{\alpha}  \langle \phi_{\alpha}(y), \phi_{\alpha}(z) \rangle$
as desired. Suppose there exists a subnet $(\phi_{\alpha_t}(y))$
such that $\Vert \phi_{\alpha_t}(y) \rVert ^2$ $\to$ $\beta$. We
need to prove that $\beta = \lVert y \rVert^2$. Clearly $\beta \leq
\lVert y \rVert ^ 2$.  If $\beta < K < \lVert y \rVert ^2$, then
there exists a $t_0$, such that, $\Vert \phi_{\alpha_t}(y) \rVert ^2
\leq K$ for all $t \geq t_0$. This implies that $\Vert
\psi_{\alpha_t} \phi_{\alpha_t}(y) \rVert ^2$ $\leq$ $\Vert
\phi_{\alpha_t}(y) \rVert ^2$ $\leq$ $K$ for all $t \geq t_0$. Since
$\psi_{\alpha_t} \phi_{\alpha_t}(y)$ $\buildrel  w^* \over \to   y$,
by Alaoglu's theorem we deduce that $\lVert y \rVert ^2$ $\leq K$,
which is a contradiction.
\end{proof}

We now generalize the notion of $W^*$-modules to the setting where
the operator algebras are $\sigma$-weakly closed algebras of
operators on a Hilbert space.  The following is the `weak* topology
variant' of the notion of {\em rigged module} studied in
\cite{DB2,BMP,DBsur,BHN} (the last paper has the most succinct
definition of these objects, and \cite{DBsur} is a survey). In
Section 4 we will prove several equivalent, but quite different
looking, characterizations of $w^*$-rigged modules.

\begin{definition} \label{wrig} Suppose that $Y$ is a dual operator
space and a right module over a dual operator algebra $M$. Suppose
that there exists a net of positive integers $(n(\alpha))$, and
$w^*$-continuous completely contractive $M$-module maps
$\phi_{\alpha} : Y \to C_{n(\alpha)}(M)$ and $\psi_{\alpha} :
C_{n(\alpha)}(M) \to Y$, with $\psi_{\alpha}( \phi_{\alpha}(y))  \to
y$ in the $w^*$-topology on $Y$, for all $y \in Y$.   Then we say
that $Y$ is a {\em right $w^*$-rigged module} over $M$.
\end{definition}

It will require some nontrivial analysis to develop the theory of
these modules from the definition given above, just as is the case
in the norm topology variant, where one needs a deep theorem of Hay
and other results \cite{BHN}. We begin by noting that an argument
similar to that in the last few lines of the proof of Lemma \ref{A},
and using basic operator space duality principles, shows that for a
$w^*$-rigged module $Y$,
\begin{equation} \label{yfor}  \Vert [y_{ij}] \Vert_{M_n(Y)}
= \sup_\alpha \; \Vert [\phi_{\alpha}(y_{ij})] \Vert , \qquad
[y_{ij}] \in M_n(Y).
\end{equation}

\begin{theorem} \label{first}  If $Y$ is a right $w^*$-rigged module over
a dual operator algebra $M$, then $w^*CB(Y)_M = {\mathcal
M}_{\ell}(Y)$ completely isometrically isomorphically, and this is a
weak*-closed subalgebra of $CB(Y)_M$. Hence $w^*CB(Y)_M$ is a dual
operator algebra, and $Y$ is a left  dual $w^*CB(Y)_M$-module.
\end{theorem}

\begin{proof}  By facts in the theory of multipliers of an operator space
(see e.g.\ \cite[Chapter 4]{BLM} or \cite{BM}), the `identity map'
is a weak* continuous completely contractive homomorphism ${\mathcal
M}_{\ell}(Y) \to CB(Y)$, which maps into $w^*CB(Y)_M$. If
$w^*CB(Y)_M$ were an operator algebra, and if $Y$ is a left operator
$w^*CB(Y)_M$-module (with the canonical action), then by the
aforementioned theory there exist a completely contractive
homomorphism $\pi : w^*CB(Y)_M \to {\mathcal M}_{\ell}(Y)$ with
$\pi(T)(y) = T(y)$ for all $y \in Y, T \in w^*CB(Y)_M$.  That is,
$\pi(T) = T$.  Thus $w^*CB(Y)_M = {\mathcal M}_{\ell}(Y)$, and it is
clear from the Krein-Smulian theorem and \cite[Theorem 4.7.4
(2)]{BLM} that $w^*CB(Y)_M$ is weak*-closed in $CB(Y)$.

We now show that $w^*CB(Y)_M$ is an operator algebra, by appealing
to the abstract characterization of operator algebras \cite[Theorem
2.3.2]{BLM}.  If $S = [S_{ij}], T = [T_{ij}] \in M_n(w^*CB(Y)_M)$,
then one may use the idea in \cite[Theorem 2.7]{DB2} or
\cite[Theorem 4.9]{BMP} to write the matrix $a = [\sum_k \, S_{ik}
T_{kj}(y_{pq})]$ as an iterated weak* limit of a product of three
matrices.  The norm of this last product is dominated by $\Vert
[S_{ij}] \Vert \Vert [T_{ij}] \Vert \Vert [y_{pq}] \Vert$.  It
follows by Alaoglu's theorem that $\Vert a \Vert \leq \Vert S \Vert
\Vert T \Vert \Vert [y_{pq}] \Vert$, and thus $\Vert S T \Vert \leq
\Vert S \Vert \Vert T \Vert$ as desired.

A similar argument shows that $Y$ is a left operator
$w^*CB(Y)_M$-module:  If $T$ is as above, and $y = [y_{ij}] \in
M_n(Y)$, then $z = [\sum_k \, T_{ik}(y_{kj})]$ may be written as a
weak* limit of a product of two matrices, the latter product having
norm $\leq \Vert T \Vert \Vert y \Vert$. Applying Alaoglu's theorem
gives $\Vert z \Vert \leq \Vert T \Vert \Vert y \Vert$, as desired.

The final assertion now follows from \cite[Lemma 4.7.5]{BLM}.
\end{proof}

\begin{theorem} \label{C}
Suppose that $Y$ is a right $w^*$-rigged module over a dual operator
algebra $M$. Suppose that $H$ is a Hilbert space, and that $\theta :
M \to B(H)$ is a  normal representation. Then $Y \otimes^{\sigma
h}_{M} H^c$ is a column Hilbert space. Moreover, the finite rank
tensors $Y \otimes H^c$ are norm dense in the latter space.
\end{theorem}

\begin{proof}
 Let
$e_\alpha = \phi_{\alpha} \psi_{\alpha}$ (notation as in Definition
\ref{wrig}).  By \cite[Lemma 2.5]{BK1} and Theorem \ref{first}, $Y
\otimes^{\sigma h}_{M} H^c$ is a left  dual ${\mathcal
M}_{\ell}(Y)$-module.  By the functoriality of the module normal
Haagerup tensor product, we obtain a net of complete contractions
$\phi_{\alpha} \otimes I_H : Y \otimes^{\sigma h}_{M} H^c \to
C_{n(\alpha)}(M) \otimes^{\sigma h}_{M} H^c$ and $\psi_{\alpha}
\otimes I_{H} : C_{n(\alpha)}(M) \otimes^{\sigma h}_{M} H^c \to Y
\otimes^{\sigma h}_{M} H^c$. Their composition $(\phi_{\alpha}
 \otimes I_H)(\psi_{\alpha} \otimes I_H) = e_\alpha
\otimes I_H$ may be regarded as the canonical left action of
$e_\alpha \in {\mathcal M}_{\ell}(Y)$ on $Y \otimes^{\sigma h}_{M}
H^c$ mentioned at the start of the proof.  Since the action is
separately weak* continuous, the composition
 converges to the
identity map on $Y \otimes^{\sigma h}_{M} H^c$ in the $w^*$-topology
(i.e.\ point-weak$^*$ on  $Y \otimes^{\sigma h}_{M} H^c)$. However,
for any $m \in \Ndb$, we have from facts in \cite{BK1} that
$$C_m(M)
\otimes^{\sigma h}_{M} H^c \cong (C_m \otimes^{\sigma h} M)
\otimes^{\sigma h}_{M} H^c  \cong C_m \otimes^{\sigma h} (M
\otimes^{\sigma h}_{M} H^c) \cong C_m \otimes^{\sigma h} H^c \cong
C_m(H^c).$$ However, $C_m(H^c)$ is a column Hilbert space. Thus by Lemma
\ref{A}, $Y \otimes^{\sigma h}_{M} H^c$ is a column Hilbert space.
The last assertion follows from \cite[Section 2]{BK1} and Mazur's
theorem that the norm closure of a convex set equals its weak
 closure.  \end{proof}

Henceforth in this section, we stick to the case that $M$ is a
$W^*$-algebra.

The following is the `Banach-module'
characterization of $W^*$-modules promised in our title. The result
may be compared with e.g.\ \cite[Corollary 8.5.25]{BLM}.

\begin{theorem}\label{B}  Let $M$ be a $W^*$-algebra.
\begin{itemize}
\item [(i)] If $Y$ is a dual Banach space and a right $M$-module,
then $Y$ is a $W^*$-module if and only if there exists a net of
integers $(n(\alpha))$, and $w^*$-continuous contractive $M$-module
maps $\phi_{\alpha} : Y \to C_{n(\alpha)}(M)$ and $\psi_{\alpha} :
C_{n(\alpha)}(M) \to Y$, with $\psi_{\alpha}( \phi_{\alpha}(y))  \to
y$ weak* in $Y$, for all $y \in Y$.
\item [(ii)]  If the conditions in {\rm (i)}
hold,
then the weak$^{*}$-limit $w^*${\rm lim}$_{\alpha} \,
\phi_{\alpha}(y)^* \phi_{\alpha}(z)$ exists in $M$ for $y$, $z$
$\in$ $Y$, and equals the $W^*$-module inner product.
\item [(iii)]  An operator
$M$-module $Y$ is $w^*$-rigged if and only if $Y$ is a $W^*$-module,
and the matrix norms for $Y$ coincide with the $W^*$-module's
canonical operator space structure.
\end{itemize}
\end{theorem}

\begin{proof}
If $Y$ is a $W^*$-module then the existence of the nets in (i) or
(iii) follow easily from e.g.\ Paschke's result \cite[Corollary
8.5.25]{BLM} or \cite[Theorem 2.1]{Bsd}.

For the other direction in (iii),
we follow the proof on p.\ 286--287 in \cite{DB1}.  
Let $\phi_{\alpha}$ and $\psi_{\alpha}$ be as in
Definition {\rm \ref{wrig}}.
We write the $k$th coordinate of
$\phi_{\alpha}$ as $x_k^{\alpha}$, where $x_k^{\alpha}$ is  a
$w^*$-continuous module map from $Y \to M$, and we write $k$th
`entry' of $\psi_{\alpha}$ as $y_k^{\alpha} \in Y$. By hypothesis we
have $\sum_{k=1}^{n(\alpha)} \, y_k^{\alpha} \, x_k^{\alpha}(y)
\buildrel w* \over \to y$ for every $y \in Y$.
Let $H$ be a Hilbert space on which $M$ is normally and faithfully
represented on, let $y, z \in Y$ and $\zeta, \eta \in H$. By Theorem
2.4 and Lemma \ref{A},
  $K = Y \otimes^{\sigma h}_{M} H$ is a column
Hilbert space.
Define two canonical maps $\Phi : Y \to B(H,K)$
and $\Psi : w^*CB_M(Y,M) \to B(K, H)$, given respectively by
$\Phi(y)(\zeta)$ = $y \otimes \zeta$ and $\Psi(f)(y \otimes
\zeta)$ = $f(y) \zeta$.  Then it is easily checked (or see
 Subsection 3.1
 for this in a more general setting), that
 $\Phi$ and $\Psi$ are weak* continuous complete
isometries.

 Let $e_{\alpha} $ = $\sum_{k=1}^{n(\alpha)} \Phi(y_k^{\alpha})
\Psi (x_k^{\alpha}) $. It is easy to check that $e_{\alpha}
\Phi(y) = \Phi(\psi_{\alpha} \phi_{\alpha}y)$ hence
$e_{\alpha} \Phi(y) \buildrel w^* \over \to \Phi(y)$ for all
$y \in Y$. Hence $e_{\alpha} (y \otimes \zeta) \to y \otimes \zeta$
weak* in $K$ for all $y \in Y, \zeta \in H$.  It follows by the last
assertion of Theorem \ref{C} that $e_{\alpha} \to I_K$ WOT in
$B(K)$.
By a similar argument to that of Theorem 4.4 in \cite{BK1}, we can
rechoose the net $(e_{\alpha})$ such that $e_{\alpha} \to I$
strongly on $K$.
Continuing to follow the proof in \cite{DB1}, one can deduce by a
small modification of the argument there, that the adjoint of any
$\Phi(y)$ $\in \Phi(Y)$ is a weak$^*$-limit of terms in
$\Psi(CB_M(Y,M))$. Thus for $z \in Y$, $\Phi(y)^* \Phi(z)$ is a
weak$^*$-limit of terms in $\Psi(w^*CB_M(Y,M)) \Phi(Y)$, and
hence is in $M$, being a weak$^*$-limit of terms in $M$.

Define $\langle y, z \rangle = \Phi(y)^* \Phi(z)$
for $y,z \in Y$. As in \cite{DB1},  $Y$ is a $C^*$-module over $M$
and the canonical $C^*$-module matrix norms coincides with the
operator space structure of $Y$, since $\Phi$ is a complete
isometry on $Y$. Since $\Phi$ is $w^*$-continuous, it follows
that the inner product on $Y$ is separately $w^*$-continuous.
Hence $Y$ is a
$W^*$-module, by e.g.\ Lemma 8.5.4 in \cite{BLM}.  This completes the proof of (iii).   

The hypothesis in (i) implies, as in the proof of (\ref{yfor}),
 that $\Vert y \Vert
= \sup_\alpha \, \Vert \phi_\alpha(y) \Vert$ for each $y \in Y$.
Then we may use (\ref{yfor}) as the definition of an operator space
structure on $Y$.  This  corresponds to an embedding of $Y$ as a
submodule of  $\oplus^\infty_\alpha \, C_{n_\alpha}(M)$, which is
easily  seen to be weak* continuous and hence a weak* homeomorphism.
Thus $Y$ becomes a dual operator module, and 
(i) and (ii) then
follow from the proof for (iii), as in \cite{DB1}, but replacing
limits by weak* limits.
\end{proof}

{\bf Remark.}   Shortly after this paper was first 
distributed, during 
conversations with Jon Kraus we found a
different proof of part of the last
result \cite{BKraus}.
We also thank Marius Junge for some answers to our
question as to whether one could find a third proof using
ultrapowers. It seems that such a proof may be more complicated than
the others.

\medskip

By a {\em weak* approximate identity} in a unital dual Banach
algebra $M$, we mean a net $\{e_t\}$ in $M$ such that $e_t $
$\buildrel w^* \over \to 1$.
A {\em weak* iterated approximate identity} for $M$
is a doubly indexed net $\{ e_{(\alpha, \beta)} \}$ (where $\beta$
and the directed set indexing $\beta$ may possibly depend on
$\alpha$), such that for each
fixed $\alpha$, the weak$^*$-limit $w^*$lim$_{\beta}$ $e_{(\alpha,
\beta)}$
exists, and $w^*{\rm lim}_{\alpha} \,
w^*{\rm lim}_{\beta} \;
 e_{(\alpha, \beta)}  \; = \; 1$.

\begin{lemma} \label{net}
A weak* iterated approximate identity for a dual Banach algebra may
be reindexed to become a weak* approximate identity.
\end{lemma}

\begin{proof}
This is a variant of the idea of the proof of \cite[Lemma 2.1]{DB2},
and we leave it to the reader.  The main new feature is that one has
to `build in' finite subsets of $M_*$ into the tuples that
constitute the elements of the new directed set.
\end{proof}

\begin{theorem} \label{okn}
Let $\{Y_i\}$ be a collection of $W^*$-modules over a $W^*$-algebra
$M$, indexed by a directed set. Let $Y$ be a dual Banach space
(resp.\ dual operator space) and a right module over $M$. Suppose
that there exist $w^*$-continuous  contractive  (resp.\ completely
contractive) $M$-module maps $\phi_{i} : Y \to Y_i$ and $\psi_{i} :
Y_i \to Y$, such that $\psi_{i}( \phi_{i}(y))  \buildrel w^{*} \over
\to y$ in $Y$, for  $y \in Y$. Then $Y$ is a $W^*$-module (resp.\ a
$W^*$-module with its canonical dual operator space structure). For
$y, z$ $\in$ $Y$, the limit $w^*${\rm lim}$_{i} \langle \phi_{i}(y),
\phi_{i}(z) \rangle$ exists in $M$ and equals the $W^*$-module inner
product $\langle y, z \rangle $.
\end{theorem}

\begin{proof}  As in Theorem \ref{B}, one can focus on the operator
space version.   For each $i$ choose nets $\phi_{\alpha _i}^i$,
$\psi_{\alpha _i}^i$ for $Y_i$ as in \ref{wrig}. Let
$\phi'_{i,\alpha _i} = \phi^i _{\alpha _i} \circ \phi_i$, and
$\psi'_{i, \alpha _i}$ = $\psi_i \circ \psi^i _{\alpha _i}$. Reindex
the net $\{\phi'_{i,\alpha _i}, \psi'_{i, \alpha _i} \}$  by Lemma
\ref{net}, so that the weak$^*$ limit of $\psi'_{i, \alpha _i}
\phi'_{i,\alpha _i}$ in $w^*CB(Y)_M$ over the new directed set
coincides with the iterated weak$^*$-limit $w^*$lim$_{i}$
$w^*$lim$_{\alpha_i}$ $\psi'_{i, \alpha _i} \phi'_{i,\alpha _i}$,
which equals $I_Y$. This gives a new asymptotic factorization of
$I_{Y}$ through spaces of form $C_n(M)$ with respect to which $Y$ is
$w^*$-rigged. Hence by Theorem \ref{B}, $Y$ is a $W^*$-module, with
the inner product
$$\langle y, z \rangle = w^*{\rm lim} \; \langle \phi^i _{\alpha _i}
(\phi_i(y)), \phi^i _{\alpha _i} ( \phi_i(z)) \rangle$$ where the
limit is taken over the new directed set.  Carefully inspecting the
directed set used in Lemma \ref{net} (a variant of the one used in
\cite[Lemma 2.1]{DB2}), it is easy to argue (and is left as an
exercise) that the last inner product equals $w^*{\rm lim}_{i} \;
\langle \phi_{i}(y), \phi_{i}(z) \rangle.$
\end{proof}

{\bf Remark.} \ The same proof as the above establishes the analogue
of the last result, but for a dual operator module $Y$  over a
unital dual operator algebra $M$, taking the $Y_i$ to be
$w^*$-rigged modules over $M$, and the $\phi_i, \psi_i$ completely
contractive (the conclusion being that $Y$ is $w^*$-rigged).

\begin{theorem} \label{P}
If $Y$ is a  right $W^*$-module over $M$, and if $Z$ is a left
(resp. right) dual operator  module over $M$, then $Y
\otimes^{\sigma h}_{M} Z$ $\cong$ $w^*CB_{M}(\overline{Y},Z) =
CB_{M}(\overline{Y},Z)$
 (resp. $Z \otimes^{\sigma h}_{M}
\overline{Y}$ $\cong$ $w^*CB(Y,Z)_{M} = CB(Y,Z)_{M}$) completely
isometrically and $w^*$-homeomorphically.
\end{theorem}

\begin{proof}  We will use facts and routine techniques from e.g.\
\cite{ER} or \cite[1.2.26, 1.6.3, Section 3.8]{BLM}.  If $T \in
B(Y,Z)_{M}$, and if $(e_{\alpha})_{\alpha \in I}$ is an orthonormal
basis for $Y$ (see \cite{Pas} or \cite[8.5.23]{BLM}), note that by
\cite[Theorem 4.2 and remark after it]{ER}
we have
$$T(y) =
T(\sum_{\alpha} \, e_{\alpha} \langle e_{\alpha} , y \rangle) =
\sum_{\alpha} \, T(e_{\alpha}) \langle e_{\alpha} , y \rangle = a
g(y), \qquad y \in Y,$$ where $a$ is the row with $\alpha$th entry
$T(e_{\alpha})$, and $g : Y \to C_I(M)$ has $\alpha$th entry the
function $\langle e_{\alpha} , \cdot \rangle$.  Thus $T$ is the
composition of `left multiplication' by  $a \in R_I(Z)$, and $g$,
both of which are weak* continuous (see e.g.\ the proof of
\cite[Corollary 8.5.25]{BLM}). Thus $w^*CB_{M}(Y,Z) = CB_{M}(Y,Z)$.
We omit the proof that $Z \otimes^{\sigma h}_{M} \overline{Y}
\cong w^*CB_M(Y,Z)$, since we prove a more general result in Theorem \ref{tpth}.
\end{proof}

\begin{corollary} \label{eqm}
In the situation of the last theorem, the tensor products
$\otimes^{\sigma h}_{M}$ occurring there coincide with Magajna's
`extended' module Haagerup tensor product $\bar{\otimes}_{hM}$ used
in {\rm \cite{Bsd}}.
\end{corollary}

It follows that in all of the results in \cite{Bsd}, all occurrences
of the `extended' module Haagerup tensor product
$\bar{\otimes}_{hM}$ may be replaced by the `normal module Haagerup
tensor product' $\otimes^{\sigma h}_{M}$. This is interesting,
since in many of these results this tensor product also coincides
with the most important and commonly used tensor product for
$W^*$-modules, the `composition' (or `fusion') tensor product $Y
\overline{\otimes}_{\theta} Z$. Thus our results gives a new way to
treat this famous `composition tensor product' (see also \cite{DH}).
Both tensor product descriptions have their own advantages:
$\bar{\otimes}_{hM}$ allows one to concretely write elements as
infinite sums of a nice form, whereas $\otimes^{\sigma h}_{M}$ has
many pleasant general properties (see \cite{EP,BK1}).

Many tensor product relations from \cite{Bsd} transfer to our
setting. For example:

\begin{corollary}  \label{fromch}
Let $Y$, $Z$ be right $W^*$-modules over $M$ and $N$
respectively,
and suppose that $\theta : M \to B(Z)$ is a normal $*$-homomorphism.
Then the `composition tensor product' $Y \overline{\otimes}_{\theta}
Z$ equals $Y \otimes^{\sigma h}_{M} Z$.
 Also, $CB(Y \overline{\otimes}_{\theta} Z )_{N}$ $\cong$
  $Y \otimes^{\sigma
h}_{M} CB(Z)_{N} \otimes^{\sigma h}_{M} \overline{Y}$ completely
isometrically and weak* homeomorphically.
\end{corollary}

\begin{proof}
The first assertion is discussed above (following
 from
 Theorem \ref{P} and \cite{Bsd}).  For the second, just as in the proof of
this result from \cite{Bsd}, Theorem \ref{P} gives
$$CB(Y \overline{\otimes}_{\theta} Z )_{N} \cong
(Y \overline{\otimes}_{\theta} Z ) \otimes^{\sigma h}_{N} (Y
\overline{\otimes}_{\theta} Z )^{-} \cong (Y \otimes^{\sigma h}_{M}
Z )  \otimes^ {\sigma h}_{N} (\overline{Z} \otimes^{\sigma h}_{M}
\overline{Y}),$$ which equals $ Y \otimes^{\sigma h}_{M}$ $ ( Z
\otimes^{\sigma h}_{N} \overline{Z} )$ $ \otimes^{\sigma h}_{M}
\overline{Y} $ $\cong$ $Y \otimes^{\sigma h}_{M} $ $B(Z)_N
\otimes^{\sigma h}_{M} $ $\overline{Y}$ (see \cite{BK1,EP}).
\end{proof}

Similarly, Theorem \ref{P} (twice) and associativity of the tensor
product, gives:

\begin{corollary}
Let $M, N$ be $W^*$-algebras, let $Y$ be a right $W^*$-module over
$M$, and let $W$ (resp.\ $Z$) be a dual operator  $N$-$M$ bimodule
 (resp.\ dual right operator  $N$-module).  Then
 $CB(Y,Z \otimes^{\sigma h}_{N} W)_M \cong Z \otimes^{\sigma h}_{N}
 CB(Y,W)_M$
 completely
isometrically and weak* homeomorphically.
\end{corollary}

\section{Some theory of $w^*$-rigged modules}

\subsection{Basic constructs}   We begin with some notation and
important constructs which will be used throughout the rest of our paper. For a
$w^*$-rigged module $Y$ over a dual operator algebra $M$, define
$\tilde{Y}$ to be $w^*CB(Y,M)_M$.   Let $\phi_{\alpha}$ and
$\psi_{\alpha}$ be as in Definition {\rm \ref{wrig}}. We write the
$k$th coordinate of $\phi_{\alpha}$ as $x_k^{\alpha}$, where
$x_k^{\alpha}$ is  a $w^*$-continuous module map from $Y \to M$, and
we write $k$th `entry' of $\psi_{\alpha}$ as $y_k^{\alpha} \in Y$.
By hypothesis we have $\sum_{k=1}^{n(\alpha)} \, y_k^{\alpha} \,
x_k^{\alpha}(y) \buildrel w^* \over \to y$ for every $y \in Y$.

 We sometimes write $\tilde{Y}$ as $X$, and denote by
 $(\cdot,\cdot)$ the canonical pairing $\tilde{Y} \times Y \to M$.
 This is completely contractive, as one may see using the idea in
 the proof of Theorem \ref{first} (the crux of the matter being
 that for $f \in \tilde{Y}, y \in Y$ we have $(f,y) = w^*{\rm lim}_\alpha \,
 \sum_{k=1}^{n(\alpha)} \, f(y_k^\alpha) x_k^\alpha(y)$, a limit
 of a product in $M$).  Let $H$ be a Hilbert space on which $M$ is normally
and faithfully (completely isometrically) represented on. Then by
Lemma \ref{C}, $K = Y \otimes^{\sigma h}_{M} H$ is a column Hilbert
space. Define two canonical maps $\Phi : Y \to B(H,K)$ and
$\Psi : \tilde{Y} \to B(K, H)$, given respectively by
$\Phi(y)(\zeta)$ = $y \otimes \zeta$ and $\Psi(f)(y \otimes
\zeta)$ = $f(y) \zeta$.  By the argument at the start of
\cite[Section 4]{BK1}, $\Phi$ is weak* continuous.  Since the
canonical map  $Y \times H^c \to K$ being completely contractive, a
routine argument gives $\Phi$ completely contractive. By the
argument on p.\ 287 in \cite{DB1}, $\Phi$ is a  complete isometry:
one obtains as in that calculation,
$$
\Vert [\phi_\alpha(y_{ij})] \Vert \leq \Vert [\Phi(y_{ij})] \Vert
,$$ so that in the limit, by (\ref{yfor}), $\Vert [ y_{ij} ] \Vert
\leq \Vert [\Phi(y_{ij})] \Vert$. The canonical weak* continuous
complete contraction
$$\tilde{Y} \otimes^{\sigma h} K^c
\cong (\tilde{Y} \otimes^{\sigma h} Y)
 \otimes^{\sigma h}_M H^c \to
M \otimes^{\sigma h}_M H^c \to H^c ,
$$ corresponds to a separately weak* continuous complete contraction
$\tilde{Y} \times K^c \to H^c$.  The map $\Psi$ is precisely the
induced weak* continuous complete contraction $\tilde{Y} \to
B(K,H)$. We leave it as an exercise that $\Psi$ is a complete
isometry.

 We define the direct sum $M \oplus^{c} Y$ as in \cite[Section
4]{BK1}.  Namely, $\theta : M \oplus Y \to B(H, K \oplus H)$ defined
by $\theta ((m,y))(\zeta) = (m \zeta , y \otimes_{M} \zeta)$,
 for $y \in Y, m \in M, \zeta \in H$, is a one-to-one
 $M$-module map,
which is a weak* continuous  complete isometry when restricted
to each of $Y$ and $M$.   We norm $M \oplus^{c} Y$  by pulling back
the operator space structure
  via $\theta$, then
 $M \oplus^{c} Y$ may be identified with the weak$^{*}$-closed
right $M$-submodule Ran$(\theta)$ of $B(H, H \oplus K)$; and hence
it is a dual operator $M$-module.

\begin{lemma} \label{sumr}  If $Y$ is a right $w^*$-rigged module
over $M$, then $M \oplus^{c} Y$ is a right  $w^*$-rigged module over
$M$.   Also, $(M \oplus^c Y) \otimes^{\sigma h}_M H^c \cong (H
\oplus K)^c$ as Hilbert spaces, for $H, K$ as above.
\end{lemma}

\begin{proof}  Define $\phi'_\alpha : M \oplus^c Y
\to C_{n(\alpha) + 1}(M)$ and $\psi'_\alpha : C_{n(\alpha) + 1}(M)
\to M \oplus^c Y$, to be $I_M \oplus \phi_\alpha$ and $I_M \oplus
\psi_\alpha$ respectively.  We also view $M \subset B(H)$, identify
$Y$ and $\Phi(Y)$, and write $n(\alpha) = n$. One may then view
$\phi'_\alpha(m,y)$, for $m \in M, y \in Y$,
 as the
matrix product of the $(n + 1) \times 2$ matrix $I_H \oplus
\Psi_{n,1}([x_k^\alpha])$ (viewed as an operator from $H \oplus K$ to 
$H^{(n+1)})$, and the $2 \times 1$ matrix with entries $m$ and
$\Phi(y)$ (viewed as an operator from $H$ to $H \oplus K$). Thus it is
clear that $\phi'_\alpha$ is completely contractive.  Similarly, we
view  $\psi'_\alpha(m,[m_k])$, for $m \in M, [m_k] \in C_n(M)$,
 as the
matrix product of the $2 \times (n + 1)$ matrix $I_H \oplus
\Phi_{1,n}([y_k^\alpha])$ (viewed as an operator from $H^{(n+1)}$ to $H
\oplus K$) and the $(n+1) \times 1$ matrix with entries $m$ and $m_k$ (viewed as an operator 
from $H$ to $H^{(n+1)}$). Thus it is clear
that $\psi'_\alpha$ is completely contractive.   It is easy to see
that $\psi'_\alpha \phi'_\alpha \to I$ weak* on $M \oplus^c Y$. So
$M \oplus^c Y$ is $w^*$-rigged.

The last assertion follows just as in \cite[Proposition 4.2]{BK1}.
\end{proof}

\begin{lemma} \label{xdos}  If $Y$ is a right $w^*$-rigged module
over $M$, then $\tilde{Y}$ is a weak* closed subspace of
$CB(Y,M)_M$.  Indeed $\tilde{Y}$ is a left $w^*$-rigged module over
$M$, which is also a dual right operator module over $w^*CB(Y)$. The
canonical map $(\cdot,\cdot) : \tilde{Y} \times Y \to M$ is
completely contractive and separately weak* continuous.
\end{lemma}

\begin{proof}  Let $P$ and $Q$ be the canonical projections
from $M \oplus^c Y$ onto $Y$ and $M$ respectively; and let $i, j$ be
the canonical inclusions of $Y$ and $M$ respectively into $M
\oplus^c Y$.  Then $\Theta(T) = j Q T i P$ defines a weak*
continuous completely contractive projection on ${\mathcal
M}_{\ell}(M \oplus^c Y) = w^*CB(M \oplus^c Y)_M$, thus its range is
weak* closed. However this range is easily seen to be completely
isometric to $w^*CB(Y,M)_M$.  Thus the latter becomes a dual
operator space, in which, from \cite[Theorem 4.7.4(2)]{BLM}, a
bounded net converges in the associated weak* topology iff it
converges point weak*.  It follows easily that $\tilde{Y}$ is a
weak* closed subspace of $CB(Y,M)_M$ (by the Krein-Smulian theorem,
or by using the fact that ${\mathcal M}_{\ell}(M \oplus^c Y)$ is
weak* closed in $CB(M \oplus^c Y)$ (see Theorem \ref{first})).

Define nets of weak* continuous maps $f \mapsto [f(y^\alpha_k)] \in
R_{n(\alpha)}(M)$, and $[m_k] \mapsto \sum_k \, m_k x^\alpha_k \in
\tilde{Y}$, then it is easy to see that with respect to these,
$\tilde{Y}$ satisfies the left module variant of Definition
\ref{wrig}.

Since $w^*CB(M \oplus^c Y)_M$ is a dual operator algebra, it is easy
to see that its `$1$-$2$-corner' $\tilde{Y}$ is a  dual right module
over its $2$-$2$-corner $w^*CB(Y)$.  We have already essentially
seen the last part.
\end{proof}

\begin{corollary} \label{cex}  We have $Y \cong w^*CB_M(\tilde{Y},M)$
completely isometrically and weak* homeomorphically, and as right $M$-modules.  That is,
$\tilde{\tilde{Y}} = Y$.  Also a bounded net $y_t \to y$ weak* in
$Y$ iff $(x,y_t) \to (x,y)$ weak* in $M$ for all $x \in \tilde{Y}$.
\end{corollary}

\begin{proof}  This is easy, using the above and the ideas
in \cite{DB2,BMP}, and routine weak* topology principles.
\end{proof}

We say that a map $T : Y \to Z$ between $w^*$-rigged modules over
$M$ is {\em adjointable} if there exists $S : \tilde{Z} \to
\tilde{Y}$ such that $(w,Ty) = (Sw,y)$ for all $y \in Y, w \in
\tilde{Z}$. The properties of adjointables in the first three
paragraphs of p.\ 389 of \cite{DB2} hold in our setting too.
Moreover by Corollary \ref{cex} and the definition of $\tilde{Y}$ we
have:

\begin{proposition}  A completely bounded module map  between $w^*$-rigged modules over $M$
is adjointable iff it is weak* continuous.
\end{proposition}

For any dual right operator $M$-modules $Y, Z$, set $\Bdb(Y,Z)$ (or $\Bdb(Y,Z)_M$) to be $w^*CB(Y,Z)_M$
and set $\Bdb(Y) = w^*CB(Y)_M$. So $\tilde{Y} = \Bdb(Y,M)$.  This use of the $\Bdb(\cdot)$ notation
may differ from that used in some of the first authors earlier papers.

We also set $N = Y \otimes^{\sigma h}_M \tilde{Y}$. Using the
canonical completely contractive and separately weak* continuous map
$(\cdot, \cdot) : \tilde{Y} \times Y \to M$, one obtains by facts in
\cite[Section 2]{BK1},  a weak* continuous completely contractive
map
$$N  \otimes^{\sigma h}
N \cong Y \otimes^{\sigma h}_M (\tilde{Y} \otimes^{\sigma h} Y)
\otimes^{\sigma h}_M \tilde{Y}
 \to Y \otimes^{\sigma h}_M M \otimes^{\sigma h}_M
 \tilde{Y} \cong N .$$
 This endows $N = Y \otimes^{\sigma h}_M \tilde{Y}$ with a
 separately  weak*
continuous completely contractive product, so that by \cite[Theorem
2.7.9]{BLM}, we have that $N$ is a dual operator algebra.  We now
show that $N$ is unital. As in \cite{BMP,DB2}, the elements
$v_\alpha = \sum_{k=1}^{n(\alpha)} \, y_k^\alpha \otimes x_k^\alpha$
are in Ball$(N)$, and for any $y \in Y, x \in \tilde{Y}$ we have in
the product above the theorem,
$$v_\alpha (y \otimes x) = \psi_\alpha(\phi_\alpha(y)) \otimes x
\, \to \, y \otimes x$$ weak* in $N$.
If $v_{\alpha_t} \to v$ is a weak* convergent subnet, then  by the
last centered formula we have $v (y \otimes x) = y \otimes x$, and
it follows that $v u = u$ for all $u \in N$.  Similarly $u v = u$.
We deduce from this that $N$ has an identity of norm $1$.  Since
such an identity is unique, we must have $v_\alpha \to 1_N$ weak*.

\begin{theorem} \label{tpth}  If $Y$ is a right
$w^*$-rigged module over $M$, and $Z$ is a right  dual operator
$M$-module, then $\Bdb(Y,Z)$ is weak* closed in $CB(Y,Z)$. Moreover,
$\Bdb(Y,Z) \cong Z \otimes^{\sigma h}_M \tilde{Y}$ completely
isometrically and weak*-homeomorphically. In particular, $\Bdb(Y)
\cong Y \otimes^{\sigma h}_M \tilde{Y}$ as dual operator algebras,
equipping the last space with the product mentioned above.
\end{theorem}

\begin{proof}
As in the second paragraph after Corollary \ref{cex}, by
\cite[Section 2]{BK1}, we have canonical weak* continuous complete
contractions
$$(Z  \otimes^{\sigma h}_M \tilde{Y}) \otimes^{\sigma h}_M Y \cong
Z \otimes^{\sigma h}_M (\tilde{Y} \otimes^{\sigma h}_M Y) \to Z
\otimes^{\sigma h}_M M \cong Z .$$ This induces a
canonical completely contractive $w^*$-continuous linear map $\theta
: Z \otimes^{\sigma h}_M \tilde{Y} \to CB(Y,Z)_M$, which satisfies
$\theta (z \otimes x) (y) $ = $ z (x,y)$, and which actually maps
into $\Bdb(Y,Z)_M$.

 In the notation above Theorem \ref{tpth}, $N = Y
\otimes^{\sigma h}_M \tilde{Y}$ is a unital dual operator algebra.
Set $W = Z \otimes^{\sigma h}_M \tilde{Y}$.
 The canonical weak* continuous maps
 $$W \otimes^{\sigma h} (Y
 \otimes^{\sigma h}_M \tilde{Y}) \cong
 Z \otimes^{\sigma h}_M (\tilde{Y} \otimes^{\sigma h} Y)
\otimes^{\sigma h}_M \tilde{Y} \to Z \otimes^{\sigma h}_M M
\otimes^{\sigma h}_M \tilde{Y} \cong W,$$ induces a separately weak*
continuous complete contraction $m : W \times N \to W$. Note that
$m(z \otimes x,1_N) = z \otimes x$ for $z \in Z, x \in \tilde{Y}$,
since $m(z \otimes x, v_\alpha) = z \otimes x \psi_\alpha
\phi_\alpha \to z \otimes x$ weak*.  Thus $m(u,1_N) = u$ for any $u
\in W$, and so $m(u,v_\alpha) \to u$ weak*.

 Now define $\mu_\alpha : CB(Y,Z)_M  \to W:
T \mapsto \sum_{k=1}^{n(\alpha)} \, T(y_k^\alpha) \otimes
x_k^\alpha$. This is a weak* continuous complete contraction. We
have $\mu_\alpha(\theta(z \otimes x)) = z \otimes x \psi_\alpha
\phi_\alpha = m(z \otimes x, v_\alpha) \to z \otimes x$ weak* for
any $z \in Z, x \in \tilde{Y}$. From the equality in the last line,
and weak* density, we have  for all $u \in W$ that
$\mu_\alpha(\theta(u)) = m(u, v_\alpha)$.  The latter, by the fact
at the end of the last paragraph, converges to $u$. Since $\Vert
\mu_\alpha(\theta(u)) \Vert \leq \Vert \theta(u) \Vert$ it follows
from Alaoglu's theorem that $\theta$ is an isometry. Similarly,
$\theta$ is a complete isometry. Since it is weak* continuous,
$\theta$ has weak* closed range, and is a weak* homeomorphism. Since
$\theta(\mu_\alpha(T)) \to T$ weak* if $T \in \Bdb(Y,Z)$, we have
now proved that Ran$(\theta) = \Bdb(Y,Z)$. Note that in the case
when $Z = Y$ we have that $\theta$ is a homomorphism, because it is
so on the weak* dense subalgebra $Y \otimes \tilde{Y}$.
\end{proof}

\begin{corollary} \label{matr} If $Y$ is a right
$w^*$-rigged module over $M$,  $Z$ is a right  dual operator
$M$-module, and $I,J$ are cardinals/sets, then \begin{itemize}
\item [(1)] $M_{I,J}(\Bdb(Y,Z)) \cong
\Bdb(C^w_J(Y),C^w_I(Z)).$
\item [(2)]  $\Bdb(C^w_I(M),Z)
\cong R_I^w(Z)$.
\item [(3)]  $\Bdb(Y,C^w_I(M)) \cong C^w_I(\tilde{Y})$.
\end{itemize}
\end{corollary}

\begin{proof}  (1) follows easily from the theorem
as on p.\ 391 in \cite{DB2}.  One should also use the fact that $\widetilde{C^w_J(Y)} =
R^w_I(\tilde{Y}) \cong \tilde{Y} \otimes^{\sigma h} R_I$, which we leave to the
reader.  (2) and (3) are immediate from (1).
\end{proof}

{\bf Remark.}   During conversations with Jon Kraus 
we were able to deduce from Theorem
\ref{tpth} that if $Y$ is a $w^*$-rigged module over
$M$ and $Z_1, Z_2$ are left dual operator $M$-modules, and if $S_t
\to S$ weak* in $w^*CB_M(Z_1,Z_2)$, then $I_Y \otimes S_t \to I_Y
\otimes S$ weak* in $CB(Y \otimes^{\sigma h}_M Z_1, Y
\otimes^{\sigma h}_M Z_2)$.   A similar result holds if $S_t
\to S$ weak* in $w^*CB(Y)_M$.  See \cite{BKraus} for details.  
This allows one to complete the `easier direction' of the analogue of one of Morita's famous
theorems: dual operator algebras
 algebras are weak* Morita equivalent iff their categories of dual operator modules are appropriately
functorially equivalent.   The `difficult direction' appears in \cite{UK1}.

\subsection{The weak linking algebra, and its representations}
\label{link}

If $Y$ is a $w^*$-rigged module over $M$, with $\tilde{Y}$, set
  $$\mathcal{L}^w =
\Big\lbrace \left[
\begin{array}{ccl}
a & x \\
y & b
\end{array}
\right] : \ a \in M, b\in \Bdb(Y), x\in \tilde{Y}, y\in Y
\Big\rbrace,  $$ with the obvious multiplication.   As in
\cite[Section 4]{BK1}, one may easily adapt the proof of the
analogous fact in \cite{BMP} to see that there is at most one possible
sensible dual operator space structure on this linking algebra.
Thus the linking
 algebra with this structure must
coincide with $\Bdb(M \oplus^c Y)$.  Another description proceeds as
follows.   Let $H$ be any Hilbert space on which $M$ is normally and
 completely isometrically represented, and set $K = Y \otimes^{\sigma
h}_M H^c$.  We saw at the start of Section 3
 the canonical maps $\Phi : Y \to B(H,K)$ and $\Psi : \tilde{Y} \to
 B(K,H)$.

\begin{proposition} \label{want}  The weak linking algebra $\mathcal{L}^w$ of the $w^*$-rigged module
may be taken to be the subalgebra of $B(H \oplus K)$ with `four
corners' $\Phi(Y), \Psi(\tilde{Y}), M$, and the weak* closure $N$ in
$B(K)$ of $\Phi(Y) \Psi(\tilde{Y})$.
In particular, $N$
is  completely isometrically
isomorphic,
via a weak* homeomorphism, to $\Bdb(Y)$.
\end{proposition}

\begin{proof}
Clearly $N$ is a weak* closed operator
algebra. Also,  $N \Phi(Y) \subset \Phi(Y)$, so that by e.g.\ 4.6.6
in \cite{BLM} we have a completely contractive homomorphism $\mu : N
\to {\mathcal M}_{\ell}(Y)$.  Conversely, since $Y$ is a dual left
operator ${\mathcal M}_{\ell}(Y)$-module by Theorem \ref{first}, so
is $K$ by \cite[Lemma 2.5]{BK1}.  Thus by the proof of \cite[Theorem
4.7.6]{BLM}, there is a normal  representation $\theta : {\mathcal
M}_{\ell}(Y) \to B(K)$.  If $y \otimes f$ denotes the obvious
operator in $w^*CB(Y)$, for $y \in Y, f \in \tilde{Y}$, then
$\theta(y \otimes f)(y' \otimes \zeta) = y f(y') \otimes \zeta =
\Phi(y) \Psi(f)(y' \otimes \zeta)$ for all $y' \in Y, \zeta \in H$.
Thus $\theta(y \otimes f) = \Phi(y) \Psi(f) \in N$. However it is
easy to see from the fact that $T \psi_\alpha \phi_\alpha \to I_Y$
weak*, that the span of such $y \otimes f$ is weak* dense in
$w^*CB(Y)_M$, and it follows that $\theta$ maps into a weak* dense
subset of $N$. Clearly $\mu(\theta(y \otimes f)) = y \otimes f$, and
so $\mu \circ \theta = I$. Thus $\theta$ is a complete isometry, and
the proof is completed by an application of the Krein-Smulian
theorem.
\end{proof}

\subsection{Tensor products of $w^*$-rigged modules} \label{ten}

 If $Y$ is a right $w^*$-rigged module over $M$, and
if $Z$ is also a right $w^*$-rigged module over a dual operator
algebra ${\mathcal R}$, and if $Z$ is a left dual operator
$M$-module, then $Y \otimes^{\sigma h}_M Z$ is a right dual operator
 ${\mathcal R}$-module (see \cite[Section 2]{BK1}).  As in the proof
of Theorem \ref{C}, we obtain a net of completely contractive right
${\mathcal R}$-module maps $\phi_{\alpha} \otimes I_Z : Y
\otimes^{\sigma h}_{M} Z \to C_{n(\alpha)}(M) \otimes^{\sigma h}_{M}
Z \cong C_{n(\alpha)}(Z)$,
 and $\psi_{\alpha} \otimes I_{Z} : C_{n(\alpha)}(Z)
\to Y \otimes^{\sigma h}_{M} Z$, such that the composition
$(\phi_{\alpha}
 \otimes I_Z)(\psi_{\alpha} \otimes I_Z) = e_\alpha
\otimes I_Z$ converges weak* to the identity map on $Y
\otimes^{\sigma h}_{M} Z$.   By the remark after Theorem \ref{okn},
$Y \otimes^{\sigma h}_{M} Z$ is a $w^*$-rigged module over
${\mathcal R}$.  In particular, if ${\mathcal R}$ is a $W^*$-algebra
then $Y \otimes^{\sigma h}_M Z$ is a $W^*$-module over ${\mathcal
R}$ by Theorem \ref{B}.

In the setting of the last paragraph (and ${\mathcal R}$ possibly
nonselfadjoint again),
\begin{equation} \label{ch1} (Y \otimes^{\sigma h}_M
Z)^{\tilde{}} = w^*CB(Y \otimes^{\sigma h}_M Z,{\mathcal
R})_{\mathcal R} \cong \tilde{Z} \otimes^{\sigma h}_M \tilde{Y} ,
\end{equation} completely isometrically and weak* homeomorphically.
We give one proof of this (another route is to use the method on p.\
402--403 in \cite{DB2}).  Note that the canonical weak* continuous
complete contractions
$$(\tilde{Z} \otimes^{\sigma h}_M \tilde{Y}) \otimes^{\sigma h} (Y \otimes^{\sigma h}_M
Z) \to \tilde{Z} \otimes^{\sigma h}_M M \otimes^{\sigma h}_M Z \to
\tilde{Z} \otimes^{\sigma h}_M Z \to {\mathcal R} ,$$ induce a weak*
continuous complete contraction $\sigma : \tilde{Z} \otimes^{\sigma
h}_M \tilde{Y} \to w^*CB(Y \otimes^{\sigma h}_M Z,{\mathcal
R})_{\mathcal R}$.  On the other hand, the complete contraction from
the operator space projective tensor product to $Y \otimes^{\sigma
h}_M Z$, induces a complete contraction $w^*CB(Y \otimes^{\sigma
h}_M Z,{\mathcal R})_{\mathcal R} \to CB(Y,CB(Z,{\mathcal R}))$ that
is easily seen to map into $CB(Y,\tilde{Z})$, and in fact into
$w^*CB(Y,\tilde{Z})_M$.  Now it is easy to check that this map
$w^*CB(Y \otimes^{\sigma h}_M Z,{\mathcal R})_{\mathcal R} \to
w^*CB(Y,\tilde{Z})_M$ is also weak* continuous.  By Theorem
\ref{tpth}, we have constructed a weak* continuous complete
contraction $\rho : w^*CB(Y \otimes^{\sigma h}_M Z,{\mathcal
R})_{\mathcal R} \to \tilde{Z} \otimes^{\sigma h}_M \tilde{Y}.$ It
is easy to check that $\rho \sigma = Id$, thus $\sigma$ is
completely isometric, and by Krein-Smulian has weak* closed range.
Any $f \in w^*CB(Y \otimes^{\sigma h}_M Z,{\mathcal R})_{\mathcal
R}$ is a weak* limit of $f \circ (\psi_\alpha \phi_\alpha \otimes
I_Z)$.  The latter function is easily checked to lie in
Ran$(\sigma)$, using the fact that for any $y \in Y$ the map $f(y
\otimes \; \cdot)$ on $Z$ is in $\tilde{Z}$.  Hence $\sigma$ has
weak* dense range, and hence is surjective, proving (\ref{ch1}).

Just as in the proof of
 Corollary \ref{fromch}, one may deduce from (\ref{ch1}) the relation
 $\Bdb(Y
\otimes^{\sigma h}_M Z) \cong Y \otimes^{\sigma h}_M
\Bdb(Z)_{\mathcal R} \otimes^{\sigma h}_M \tilde{Y}$.

In fact the weak* variants of all the theorems in Section 6 of
\cite{DB2} are valid (and we will perhaps present these elsewhere).
 In next subsection, we merely focus on 6.8 from that
paper, which we shall need at the end of the next section.

 \subsection{The $W^*$-dilation} \label{wdil}

This important tool is
 a canonical `$W^*$-module `envelope' of a $w^*$-rigged module $Y$
 over $M$.
  If ${\mathcal R}$ is a $W^*$-algebra containing $M$
as a weak*-closed subalgebra with $1_{\mathcal R} = 1_M$, then $E =
Y \otimes^{\sigma h}_M {\mathcal R}$ is a $W^*$-module over
${\mathcal R}$ by \ref{ten}, and it is called a {\em $W^*$-dilation} of $Y$.  We may
identify $Y$ with $Y \otimes 1$. This is a completely isometric
weak* homeomorphic identification, since by (\ref{yfor}) we have for
$[y_{ij}] \in M_n(Y)$ that
$$\Vert [ y_{ij} \otimes 1 ] \Vert_{M_n(E)} = \sup_\alpha \, \Vert [
(\phi_{\alpha} \otimes I_{Z})(y_{ij} \otimes 1) ] \Vert =
 \sup_\alpha \, \Vert [ \phi_{\alpha}(y_{ij}) ] \Vert =
 \Vert [ y_{ij} ] \Vert_{M_n(Y)} .$$
Thus every $w^*$-rigged module over $M$ is a weak* closed
$M$-submodule of a $W^*$-module over ${\mathcal R}$.  Usually we
assume that ${\mathcal R}$ is generated as a $W^*$-algebra by $M$.

Similarly, it is easy to see that  $\tilde{Y}$ is a weak* closed
left $M$-submodule of ${\mathcal R} \otimes^{\sigma h}_M \tilde{Y}$.
By (\ref{ch1})  above, ${\mathcal R} \otimes^{\sigma h}_M \tilde{Y}
= \tilde{E}$, which in turn is just the conjugate $C^*$-module (see
e.g.\ \cite[8.1.1 and and 8.2.3(2)]{BLM}) $\bar{E}$ of $E$. We claim
that $\Bdb(Y)$ may be regarded as  a weak* closed subalgebra of
$\Bdb(E)$ having a common identity element.  By a principle we have
met several times now (e.g.\ in the proof of Proposition \ref{want}),
there is a canonical weak* continuous completely contractive unital
homomorphism $\Bdb(Y) \to \Bdb(E)$. However, since $Y \cong Y
\otimes 1 \subset E$ as above, it is easy to see that the last
homomorphism  is a completely isometric weak* homeomorphism.  Thus
we have established the variant in our setting of \cite[Theorem
6.8]{DB2}.

The $W^*$-dilation is studied in a more general setting in
\cite{UK1,UK}.

\subsection{Direct sums} \label{dirs}

If $Y$ is a $w^*$-rigged module over $M$, and if $P \in \Bdb(Y)$ is
a contractive idempotent, then it is easy to see from the Remark
after Theorem \ref{okn}, that $P(Y)$ is a $w^*$-rigged module over
$M$, called an {\em orthogonally complemented submodule} of $Y$.

As in the discussion at the top of p.\ 409 of \cite{DB2}, if $\{ Y_k
\}_{k \in I}$ is a collection of $w^*$-rigged modules over $M$, and
if $E_k = Y_k \otimes^{\sigma h}_M {\mathcal R}$ is the
$W^*$-dilation of $Y_k$, for a $W^*$-algebra ${\mathcal R}$
containing $M$, then we can define the `column direct sum'
$\oplus^c_{k \in I} \; Y_k$, to be $\oplus^c_{k \in I} \; Y_k = \{
(y_k) \in \oplus^c_{k \in I} \; E_k : y_k \in Y_k \; \textrm{ for
all} \;  k \in I \}$, where $\oplus^c_{k \in I} \; E_k$ is the
$W^*$-module direct sum (see \cite{Pas} or \cite[8.5.26]{BLM}). A
key principle, which is used all the time when working with this direct sum,
 is the following.
Consider the directed set of finite subsets $\Delta$ of $I$, and for
$z \in \oplus^c_{k \in I} \; E_k$, write $z_\Delta$ for the tuple
$z$, but with entries $z_k$ switched to zero if $k \notin \Delta$.
Then $(z_\Delta)_\Delta$ is a net indexed by $\Delta$, which
converges weak* to $z$.  For example, it follows from this principle
that $\oplus^c_{k \in I} \; Y_k$ is the weak* closure inside
$\oplus^c_{k \in I} \; E_k$ of the `finitely supported tuples'
$(y_k)$ with $y_k \in Y_k$ for all $k$.

\begin{theorem}
\label{L22}  If $\{ Y_k : k \in I \}$ is a collection of
$w^*$-rigged modules over $M$, then  $\oplus^c_{k \in I} \; Y_k$ is
again a $w^*$-rigged module over $M$.  Also $\widetilde{\oplus^c_{k \in I} \; Y_k}
= \oplus^r_{k \in I} \; \tilde{Y_k}$.
\end{theorem}

\begin{proof}  We first observe that this is easy if $I$ is finite.  The reader
may just want to consider the case of two modules, which is similar
to the proof of Lemma \ref{sumr}. In fact one may use Definition
\ref{5} below to see quickly that $Y_1 \oplus^c Y_2$ is $w^*$-rigged
if $Y_1$ and $Y_2$ are: note that if $E = E_1 \oplus^c E_2$ in the
notation of the last paragraph, and if $(z_1 , z_2) \in E$ is such
that $\langle (z_1 , z_2) | (y_1,y_2) \rangle = \langle z_1 | y_1
\rangle + \langle z_2 | y_2 \rangle \in M$ for all $y_1,y_2 \in Y$,
then $z_k^* \in X_k = \tilde{Y_k}$ by  \ref{5}. It is then easy to see that the
conditions of Definition \ref{5} are satisfied, so that $Y_1 \oplus^c Y_2$
is $w^*$-rigged.

Suppose that $I$ is infinite.  In the notation of the paragraph
above the theorem, we have that $\oplus^c_{k \in \Delta} \; Y_k$ is
$w^*$-rigged by the last paragraph.  There are canonical maps
$\phi_\Delta$ and $\psi_\Delta$ between $Y = \oplus^c_{k \in I} \;
Y_k$ and $\oplus^c_{k \in \Delta} \; Y_k$.  Namely, $\phi_\Delta$ is
essentially the map $z \mapsto z_\Delta$, and $\psi_\Delta$ is the
`inclusion', indeed $\phi_\Delta \circ \psi_\Delta = Id$.  It is
easy to see that these maps are completely contractive and weak*
continuous, since when one tensors them with $I_{\mathcal R}$ they
have these properties. Also, $\psi_\Delta \circ \phi_\Delta \to I_Y$
weak*, using the principle above the theorem that $z_\Delta \to z$.
It follows from the Remark after Theorem \ref{okn}, that $Y$ is
$w^*$-rigged. We leave the last relation to the reader.
  \end{proof}

The following universal property shows that the direct sum
$\oplus^c_{k \in I} \; Y_k$ does not dependent on the specific
construction of it above:

\begin{theorem}
\label{L2} Suppose that $\{ Y_k \}_{k \in I}$ is a collection of
dual operator modules over $M$, that $Y$ is a fixed $w^*$-rigged
module over $M$, and that there exist weak* continuous completely contractive
$M$-module maps $i_k : Y_k \rightarrow Y$, $\pi_k : Y \rightarrow
Y_k$ with $\pi_k \circ i_m = \delta_{km} Id_{Y_m}$ for all $k,m$.
Here $\delta_{km}$ is the Kronecker delta.
Then each $Y_k$ is $w^*$-rigged, and $Y$ is completely isometrically
weak* homeomorphically, $M$-isomorphic  to the column direct sum $Z
\oplus^c (\oplus^c_k Y_k)$ defined above, where $Z$ is a submodule
of $Y$ which is also $w^*$-rigged. If $\sum_{k \in I} \, i_k \pi_k =
I_Y$ in the weak*-topology of $\Bdb(Y)$, then $Z = (0)$.
\end{theorem}

\begin{proof}
The ranges $i_k(Y_k)$ are orthogonally complemented submodules of
$Y$, and hence they are $w^*$-rigged, and so is $Y_k$. The partial
sums of  $R = \sum_k \, i_k \pi_k$  form an increasing net of
contractive projections in the dual operator algebra $\Bdb(Y)$,
indexed by the finite subsets of $I$ directed upwards by inclusion.
Hence it converges in the weak* topology in $\Bdb(Y)$ to a
contractive projection $R \in \Bdb(Y)$. Let $Z = {\rm Ran}(I-R)$,
which again is $w^*$-rigged. Define $Z \oplus^c (\oplus^c_k \; Y_k)$
as above the theorem, a weak* closed $M$-submodule of the
$W^*$-module direct sum $(Z \otimes^{\sigma h}_M {\mathcal R})
\oplus^c (\oplus^c_k \; (Y_k \otimes^{\sigma h}_M {\mathcal R}))$.
Tensoring all maps with $I_{\mathcal R}$, we obtain maps back and
forth between $Y_k \otimes^{\sigma h}_M {\mathcal R}$ and $Y
\otimes^{\sigma h}_M {\mathcal R}$, and between $Z \otimes^{\sigma
h}_M {\mathcal R}$ and $Y \otimes^{\sigma h}_M {\mathcal R}$,
satisfying the hypotheses of \cite[Theorem 2.2]{Bsd}.   Note that
$i_k \pi_k \in \Bdb(Y)_M$, and $Y \otimes^{\sigma h}_M {\mathcal R}$
is a left dual operator $\Bdb(Y)_M$-module (since $Y$ is).  It
follows that $\sum_k \, (i_k \pi_k \otimes I_{\mathcal R}) = R
\otimes I_{\mathcal R}$, so that $\sum_k \, (i_k \pi_k \otimes
I_{\mathcal R}) + (I- R) \otimes I_{\mathcal R} = I$.  From
\cite[Theorem 2.2]{Bsd}, it follows that the canonical map is a
completely isometric weak* homeomorphic, ${\mathcal R}$-isomorphism
between $(Z \otimes^{\sigma h}_M {\mathcal R}) \oplus^c (\oplus^c_k
\; (Y_k \otimes^{\sigma h}_M {\mathcal R}))$ and $Y \otimes^{\sigma
h}_M {\mathcal R}$. Its restriction  to  the copy of $Z \oplus
(\oplus^w_n Y_n)$ is the desired map.
\end{proof}

As in \cite[Section 7]{DB2} it follows that the column direct sum is
associative and commutative.   We also have the obvious variant of
\cite[Theorem 7.4]{DB2} valid in our setting, concerning the direct
sum $\oplus_k \, T_k$ of maps  $T_k \in \Bdb(Y_k,Z_k)$.  Again, the
proof of this is now familiar: apply the $W^*$-module case of this
result to the maps $T_k \otimes I_{\mathcal R}$ between the
$W^*$-dilations, and then restrict to the appropriate subspace.
Also, we obtain from Theorem \ref{L2} and functoriality of the
tensor product $\otimes^{\sigma h}_M$, as in \cite[p.\ 411]{DB2},
both left and right distributivity of this tensor product
$\otimes^{\sigma h}_M$ over column direct sums of $w^*$-rigged
modules: $$(\oplus^c_k \; Y_k) \otimes^{\sigma h}_M Z \cong
\oplus^c_k \; (Y_k \otimes^{\sigma h}_M Z) ,$$ and
$$Y \otimes^{\sigma h}_M (\oplus^c_k \; Z_k)
\cong  \oplus^c_k \; (Y \otimes^{\sigma h}_M Z_k) .$$
All spaces in these formulae are right $w^*$-rigged modules, and $Z$
and $Z_k$ are also left dual operator $M$-modules. For the last
formula, one may use
the fact that
if $T_t \to T$ weak* in $\Bdb(Z)$, then $I_Y \otimes T_{t} \to I
\otimes T$ weak* in $\Bdb(Y \otimes^{\sigma h}_M Z)$. Indeed, if we
have a weak* convergent subnet $I_Y \otimes T_{t_\mu} \to S \in
\Bdb(Y \otimes^{\sigma h}_M Z)$, then $S(y \otimes z) = y \otimes
T(z)$ for $y \in Y, z \in Z$. Since finite rank tensors are weak*
dense, we have $S = I \otimes T$, and it follows that $I_Y \otimes
T_{t} \to I \otimes T$ weak*.
Full details are left to the reader.

 \medskip

 {\bf Remark.}  Theorem \ref{L2} also shows that the definition
 from Section 3 of  $M \oplus^c Y$, agrees with the column
 direct sum in the present subsection.  Thus the last relation in
Lemma \ref{sumr} is a simple special case of the second last
centered formula.

\section{Equivalent definitions of $w^*$-rigged modules}

\subsection{}  The reader may prefer some of the following
four descriptions of $w^*$-rigged modules, each of which involve a
{\em pair} $X, Y$ of modules. In each case, the first paragraph of
the subsection  constitutes the alternative definition.  In the
second and following paragraphs the equivalence with the original
definition is sketched.  One must show that every $w^*$-rigged
module $Y$ satisfies (or is completely isometrically, weak*
homeomorphically, $M$-isomorphic to a module which satisfies) the
given alternative description; and that conversely any $Y$
satisfying the description is $w^*$-rigged, and that moreover $X
\cong \tilde{Y}$.

We will be a little informal in this section, since the objectives
here are quite clear---we are just adapting four theorems from
\cite[Section 5]{DB2} to the weak* topology setting of the present
paper.  The reader will  easily be able to
 add any missing detail.

\subsection{Second definition of a $w^*$-rigged
module} \label{2}

 Fix two unital dual operator algebras $M$ and $N$, and two
dual operator bimodules $X$ and $Y$, with $X$ an $M$-$N$-bimodule
and $Y$ an $N$-$M$-bimodule.  We also assume that there exists a
separately weak$^{*}$-continuous completely contractive $M$-bimodule
map $( \cdot,\cdot ) : X \times Y \to M$ which is balanced over $N$,
and a
 separately weak$^{*}$-continuous completely contractive  $N$-bimodule
  map  $\lbrack \cdot, \cdot \rbrack : Y \times X \to N$ which is balanced over $M$,
such that $(x,y)x' = x [y,x']$ and $y' (x,y) = [y',x] y$ for $x,x'
\in X, y,y' \in Y$; and such that $\lbrack \cdot, \cdot \rbrack$
induces a weak* continuous quotient map $Y \otimes^{\sigma h} X \to
N$.

As in \cite[Section 5]{DB2}, any $w^*$-rigged module in the sense of
Definition \ref{wrig}, satisfies the conditions in the last
paragraph, with $N = \Bdb(Y)$ (or $N = Y \otimes^{\sigma h}_M
\tilde{Y}$), and $X = \tilde{Y}$, by our earlier results.
Conversely, given the conditions in the last paragraph, suppose that
$u \in {\rm Ball}(Y \otimes^{\sigma h} X)$ maps to $1_N$, and that
$(f_s)$ is a net of finite rank tensors in Ball$(Y \otimes_{\rm h}
X)$ which converges weak* to $u$ (using \cite[Corollary 2.6]{BK1}).
The image of $f_s$ in $N$ converges weak* to $1_N$. From this we see
that $Y$ satisfies Definition \ref{wrig} (as in similar assertions
in \cite{DB2} (see e.g.\ bottom of p.\ 384 there)). Moreover, a by
now routine modification of the last two paragraphs of the proof of
\cite[Theorem 4.1]{BMP}, one sees that the canonical map $X \to
w^*CB(Y,M)_M$ is a weak* continuous surjective complete isometry.
That is, $X \cong \tilde{Y}$ as dual operator $M$-modules. We have a
canonical weak* continuous complete quotient map $\theta : Y
\otimes^{\sigma h}_M \tilde{Y}\to N$. A simple modification of the
last paragraph of the proof of Theorem \ref{tpth}, which is
essentially the proof of ($\Leftarrow$) in \cite[Theorem 3.3]{BK1},
shows that $\theta$ is a complete isometry, so that $[\cdot , \cdot
]$ induces a weak* homeomorphic complete isometry $Y \otimes^{\sigma
h}_M \tilde{Y} \cong N$.

\subsection{Third definition of a $w^*$-rigged
module} \label{3} A pair consisting of a dual left $M$-module $X$,
and a dual right $M$-module $Y$, with a separately weak* continuous
completely contractive pairing $(\cdot , \cdot) : X \times Y \to M$,
such that if we equip $N = Y \otimes^{\sigma h}_M X$ with the
canonical
 separately  weak*
continuous completely contractive product induced by $(\cdot ,
\cdot)$, as in the discussion above Theorem \ref{tpth}, then this
(dual operator) algebra has an identity of norm $1$.  We also assume
that the canonical actions of $N$ on  $Y$ and on $X$ are
non-degenerate (that is, $1_N y = y, x 1_N = x$ for $y \in Y, x \in
X$).

Again, clearly any $w^*$-rigged module in the earlier sense,
satisfies the conditions in the last paragraph, by Theorem
\ref{tpth} and the remarks above it. Conversely, suppose that $X, Y,
(\cdot , \cdot)$ are as in the last paragraph; we shall verify the
conditions of Definition \ref{2}.  It is by-now-routine to see that
$X, Y$ are dual operator modules over $N$. To see that $(\cdot ,
\cdot)$ is $N$-balanced, one shows that for $x \in X, y \in
Y$, the two weak* continuous functions $(x, \, \cdot \; y)$ and $(x
\; \cdot \, , y)$ on $N$, are equal on the weak* dense subset $Y
\otimes X$ of $N$. The rest is obvious.

\subsection{Fourth description of $w^*$-rigged
modules} \label{4} Let $M, N$ be weak* closed unital subalgebras
 of $B(H)$ and $B(K)$ respectively, for Hilbert spaces $H, K$, and
let $X \subset B(K,H), Y \subset B(H,K)$ be weak* closed subspaces,
such that the associated subset ${\mathcal L}$ of $B(H \oplus K)$ is
a subalgebra of $B(H \oplus K)$, for Hilbert spaces $H, K$.  This is
the same as specifying a list of obvious algebraic conditions, such
as $X Y \subset M$.  Assume in addition that
the weak* closure $N$ in $B(K)$ of $Y X$, possesses a net $(e_t)$
with terms of the form $y x$, for $x \in {\rm Ball}(C_n(X))$ and $y
\in {\rm Ball}(R_n(Y))$, such that $e_t \to 1_N$ weak*.

That every $w^*$-rigged module $Y$ is essentially of this form,
follows by replacing $Y$ and $\tilde{Y}$ by $\Phi(Y)$ and $X =
\Psi(\tilde{Y})$ respectively, and looking at the weak linking algebra ${\mathcal
L}^w$ at the end of Section 3.2.
Let $e_{\alpha}  = \sum_{k=1}^{n(\alpha)} \Phi(y_k^{\alpha})
\Psi (x_k^{\alpha}) $. It is easy to check that $e_{\alpha}
\Phi(y) = \Phi(\psi_{\alpha} \phi_{\alpha}y)$, hence
$e_{\alpha} \Phi(y) \buildrel w^* \over \to \Phi(y)$ for all
$y \in Y$. Hence $e_{\alpha} (y \otimes \zeta) \to y \otimes \zeta$
weak* in $K$ for all $y \in Y, \zeta \in H$.  It follows by the last
assertion of Theorem \ref{C} that $e_{\alpha} \to I_K$ WOT in
$B(K)$.
This gives the condition in the last paragraph. Conversely,
given the setup in the last paragraph, we will verify the conditions
of Definition \ref{2}.  The canonical map $\theta : Y
\otimes^{\sigma h} X \to N$ is completely contractive and weak*
continuous, we need to show it is a quotient map.  If $T \in {\rm
Ball}(N)$, and if we write the $x$ and $y$ in the last paragraph as
$x = [x_k], y = [y_k]$, then $u_t = \sum_k \, T y_k \otimes x_k \in
{\rm Ball}(Y \otimes^{\sigma h} X)$.  Consider a weak* convergent
subnet $u_{t_\beta } \to u \in {\rm Ball}(Y \otimes^{\sigma h} X)$.
Then $\theta(u_{t_\beta }) \to \theta(u)$. On the other hand,
$\theta(u_{t_\beta}) = T e_{t_\beta} \to T$ weak*.  So $T =
\theta(u)$, so that $\theta$ is a quotient map.

\subsection{Fifth definition of a $w^*$-rigged
module}  \label{5} Let ${\mathcal R}$ be a $W^*$-algebra containing
$M$ as a weak* closed subalgebra with $1_{\mathcal R} = 1_M$, and
suppose that $Z$ is a right $W^*$-module over ${\mathcal R}$, and
that $Y$ is a weak* closed $M$-submodule of $Z$.  Define $W = \{ z
\in Z : \langle z \vert y \rangle \in M \}$, and set $N$ to be the
weak* closure in $\Bdb(Z)_{\mathcal R}$ of the span of terms of the
form $| y \rangle \langle w |$ for $y \in Y, w \in W$ (here $| y
\rangle \langle w |$ is the obvious `rank one' operator $z \mapsto y
\langle w , z \rangle$ on $Z$. Suppose that there is a net $(e_t)$
converging to $I_Z$ weak* in $\Bdb(Z)$, with terms of the form $e_t
= \sum_{k=1}^n \, | y_k \rangle \langle w_k |$, where $y_k \in Y,
w_k \in W$ with $\sum_{k=1}^n \, | y_k \rangle \langle y_k | \leq 1$
and $\sum_{k=1}^n \, | w_k \rangle \langle w_k | \leq 1$.

We claim that under the hypotheses in the last paragraph, $Y$  is
$w^*$-rigged, $\tilde{Y} = \{ \bar{w} \in \bar{Z} : w \in W \}$, and
$\Bdb(Y)_M \cong N$.   To see this, we follow the proof of
\cite[Theorem 5.10]{DB2}, working inside the linking $W^*$-algebra
${\mathcal L}^w(Z)$ for $Z$, where all inner products and module
actions become concrete operator multiplication. Note first that $W$
is a weak* closed right $M^*$-submodule of $Z$, and hence  $X = W^*$
is a weak* closed left $M$-submodule of $Z^*$. The subspace of
${\mathcal L}^w(Z)$ with four corners $M, X, Y, N$, is a weak*
closed subalgebra, and it is easy to see that the criteria of the
Definition \ref{4} are met, for these subspaces of ${\mathcal
L}^w(Z)$. Hence the criteria of Definition \ref{2}  are met, and we
are done by facts from that place.

 Conversely, to see that every $w^*$-rigged module $Y$ is
  essentially of this form,
 set  $Z = Y \otimes^{\sigma h}_M {\mathcal R}$, which we saw in
 \ref{wdil}
 is a $W^*$-module
over ${\mathcal R}$, containing $Y$ as a weak* closed $M$-submodule.
Also we saw  that $\tilde{Y} \subset \bar{Z}$ (resp.\
$\Bdb(Y)_M \subset \Bdb(Z)_{\mathcal R}$) as a weak* closed
$M$-submodule (resp.\ weak* closed subalgebra with common
identity).  Now apply a simple variant of the argument in the last
paragraph of \cite[p.\ 405]{DB2}.

\section{Examples}

\begin{itemize}
\item [(1)]  As we saw in Section 2, $W^*$-modules  are $w^*$-rigged.
Thus so are {\em  WTROs}, where a WTRO is  a weak* closed space $Z$
of Hilbert space operators  with $Z Z^* Z \subset Z$ (see
\cite[8.5.11 and 8.5.18]{BLM}).
\item [(2)]  For finite dimensional  modules over a finite
dimensional operator algebra $M$, the notions of `rigged' and
`$w^*$-rigged' coincide, as is easily seen from Definition
\ref{wrig}.
\item [(3)]   By e.g.\ \ref{2} and \cite[Theorem 3.3]{BK1}, every weak* Morita
equivalence bimodule in the sense of \cite{BK1} is $w^*$-rigged. In
Section 3 of that paper, a long list of examples of these bimodules
is given.  Indeed a weak* Morita equivalence bimodule is essentially
the same thing as a `left-right symmetric' variant of Definition
\ref{2} (that is, we also assume there that  $( \cdot, \cdot )$
induces a weak* continuous quotient map $X \otimes^{\sigma h} Y \to
M$).

There are simple examples of $w^*$-rigged modules  which give rise
to no kind of weak* Morita equivalence (in contrast to the
$W^*$-module case). For example, consider $Y = R_2$, a right
$w^*$-rigged module over the upper triangular $2 \times 2$ matrices.
A partial result in the positive direction here: if $Y$ is a
$w^*$-rigged $M$-module which is {\em $w^*$-full}, that is the span
of the range of $(\cdot,\cdot)$ is weak* dense in $M$, and if
${\mathcal R}$ is a $W^*$-algebra generated by $M$, then the
$W^*$-dilation $E = Y \otimes^{\sigma h}_M {\mathcal R}$ gives a von
Neumann algebraic Morita equivalence (see \cite{Rief} or
\cite[Section 8.7]{BLM}) between ${\mathcal R}$ and $\Bdb(E)$.  This
will follow if $E$ is $w^*$-full over ${\mathcal R}$ (see
\cite[8.5.12]{BLM}). To this end, note $\bar{E} = \tilde{E} =
{\mathcal R} \otimes^{\sigma h}_M \tilde{Y}$ by (\ref{ch1}).  Thus
$\tilde{Y} Y$, and therefore also $M$, is contained in the weak*
closure of $\tilde{E} E$. So $E$ is $w^*$-full, since the latter is
an ideal of ${\mathcal R}$, and because $M$ generates ${\mathcal
R}$.
\item [(4)]  The second dual of a rigged module over an operator algebra
$A$ is $w^*$-rigged over $A^{**}$.  This is evident by taking the
second dual of all objects in the definition of a rigged module from
\cite{BHN} say (note that $C_n(A)^{**} = C_n(A^{**})$ by basic
operator space duality).
\item [(5)]  If $P$ is a weak*-continuous completely contractive
idempotent $M$-module map on $C^w_I(M)$, for a cardinal/set $I$,
then Ran$(P)$ is a $w^*$-rigged module (see the first paragraph of
Section  \ref{dirs}).
\item [(6)] Examples of $w^*$-rigged may be built analogously to the
rigged modules constructed in \cite{BJ} (see e.g.\ the end of Section 6
 there).
\item [(7)]   In \cite{BK1} (see also \cite{UK})
we defined a stronger
variant of weak* Morita equivalence (see Example (3) above)
called {\em weak Morita equivalence}.
In just the same way, we can define
a subclass of $w^*$-rigged modules:
we say that an $M$-module is {\em weakly rigged}
if it satisfies
Definition 3.2 of \cite{BK1}, but
with the phrase `(strong) Morita context' replaced by `(P)-context' (see
\cite[p.\ 20]{BMP}).  An adaption of the proof of \cite[Corollary
3.4]{BK1} shows that weakly rigged modules are $w^*$-rigged.  One
may then show  that any weakly rigged module pair $(Y,X)$ arises as
a weak* closure of a rigged module situation, just as in Example (2)
after Definition 3.2 of \cite{BK1}, but dropping the requirement on
the cai for $A$ there. This proceeds by showing that the linking
algebra for the `subcontext' is a weak* dense subalgebra of the weak
linking algebra for $(Y,X)$ (this uses 6.10 in \cite{DB2}, see the
argument above 4.1 in \cite{BK1}).
\item [(8)]  Let $Z$ be any WTRO (see (1)), and suppose that
$Z^* Z$ is contained in a dual operator algebra $M$. Then $Y =
\overline{ZM}^{w^*}$ is a $w^*$-rigged $M$-module.  We call this
example a {\em WTRO-rigged module}.  We also have $\tilde{Y} \cong
\overline{M Z^*}^{w^*}$.  To see all this, denote the last space by
$X$, and set $N$ to be the weak* closure of $Z M Z^*$, a dual
operator algebra containing $Z Z^*$. If $(e_t)$ is the usual
approximate identity for $Z Z^*$ with terms of the form
$\sum_{k=1}^n \, z_k z_k^*$, then it is routine to see that $(e_t)$
converges weak* to an identity $1_N$ for $N$. Now it is easy to
check that $M, Y, X, N$ satisfy Definition \ref{4}, and we are done.

 We remark that the above is a generalization of Eleftherakis'
 recent notion of  {\em TRO-equivalence} (see e.g.\ \cite{El,EP}).
  Indeed,
 a WTRO-rigged module gives a TRO-equivalence of $M$ and $N$ iff
 the identity $e$ of the weak* closure of $Z^* Z$ is $1_M$.
 For the most difficult part of this, note that if the latter
 holds, and if $f_s \to e$ weak* with $f_s \in Z^* Z$, then
 any $m \in M$ is an iterated  weak* limit of the $f_s m f_{s'}$,
 and it follows that $M$ equals the weak* closure of $Z^* N Z$.
\item [(9)]  The {\em selfdual} rigged modules over a
dual operator algebra $M$, considered at the start of the last
section in \cite{BM}, together with their unique dual space
structure making $(\cdot,\cdot)$  separately weak* continuous (see
\cite[Lemma 5.1]{BM}),
 are $w^*$-rigged.  Indeed it is easy to see from the last mentioned
 continuity that Definition \ref{wrig} is satisfied.
  \end{itemize}

In \cite{BKraus} we prove that classes {\rm (5)} and {\rm (8)} above
coincide, and also equal the  $w^*$-rigged module 
 direct sums of modules of the form $p_i M$, for projections $p_i \in M$.  
This is the analogue of a famous theorem due to Paschke
(see \cite{Pas} or \cite[Corollary 8.5.25]{BLM}).  We call these 
{\em projectively $w^*$-rigged modules}.  
We show that they constitute a proper subclass of the $w^*$-rigged modules, and that 
they have some strong properties.  For example, if $Y$ is a weak* Morita equivalence $M$-$N$-bimodule,
over dual operator algebras
$M$ and $N$, which is both left and right projectively
$w^*$-rigged, then $M$ and $N$ are stably isomorphic (that is,
$M_I(M) \cong M_I(N)$ completely isometrically for some cardinal $I$).
This is a reprise of the main result of \cite{EP}.

\end{document}